\documentstyle{amsppt}
\magnification1200
\pagewidth{6.5 true in}
\pageheight{9.25 true in}
\NoBlackBoxes
\topmatter
\title 
The fourth moment of Dirichlet $L$-functions  
\endtitle
\author 
K. Soundararajan
\endauthor 
\address 
Department of Mathematics, University of Michigan, Ann Arbor, MI 48109, USA
\endaddress
\email 
ksound\@ math.lsa.umich.edu
\endemail

\endtopmatter

\def\sumstar{\sideset \and^{*} \to \sum}

\def\phi{\varphi}

\def\Lam{\Lambda}
\def\cbar{\overline{\chi}}

\document

\head 1. Introduction  \endhead

\noindent In [2], D.R. Heath-Brown showed that 
$$
\sumstar_{\chi\pmod q} |L(\tfrac 12,\chi)|^4 
= \frac{\phi^*(q)}{2\pi^2} \prod_{p|q} \frac{(1-p^{-1})^3}{(1+p^{-1})} 
(\log q)^4 + O(2^{\omega(q)} q(\log q)^{3}). \tag{1.1}
$$
Here $\sumstar$ denotes summation over primitive characters 
$\chi \pmod q$, $\phi^*(q)$ denotes the number of 
primitive characters $\pmod q$, and $\omega(q)$ 
denotes the number of distinct prime factors of $q$.  
Note that $\phi^*(q)$ is a multiplicative function 
given by $\phi^*(p) =p-2$ for primes $p$, and 
$\phi^*(p^k) = p^k (1-1/p)^2$ for $k\ge 2$ (see Lemma 1 below).  
Also note that when $q\equiv 2\pmod 4$ there are no primitive 
characters $\pmod q$, and so below we will assume 
that $q\not\equiv 2 \pmod 4$.  For $q\not\equiv 2 \pmod 4$ it is 
useful to keep in mind that the main term in (1.1) is 
$\asymp q (\phi(q)/q)^6  (\log q)^4$. 

Heath-Brown's result represents a $q$-analog of Ingham's fourth moment 
for $\zeta(s)$: 
$$
\int_0^T |\zeta(\tfrac 12+it)|^4 dt \sim \frac{T}{2\pi^2}(\log T)^4.
$$
When $\omega(q) \le (1/\log 2 -\epsilon) \log \log q$ (which holds 
for almost all $q$) the error term in (1.1) is dominated by the 
main term and (1.1) gives the $q$-analog of Ingham's 
result.  However if $q$ is even a little 
more than `ordinarily composite', with 
$\omega(q) \ge (\log \log q)/\log 2$, then the error term 
in (1.1) dominates the main term.  In this note we remedy this, 
and obtain an asymptotic formula valid for all large $q$. 

\proclaim{Theorem}  For all large $q$ we have 
$$
\sumstar_{\chi \pmod q} |L(\tfrac 12,\chi)|^4 
= \frac{\phi^*(q)}{2\pi^2} \prod_{p|q} \frac{(1-p^{-1})^3}{(1+p^{-1})} 
(\log q)^4 \Big( 1+ 
O\Big(\frac{\omega(q)}{\log q} \sqrt{\frac{q}{\phi(q)}} \Big) \Big)
+ O(q (\log q)^{\frac 72}).
$$
\endproclaim

Since $\omega(q) \ll \log q/\log \log q$, and $q/\phi(q) \ll 
\log \log q$, 
we see that $(\omega(q)/\log q) \sqrt{q/\phi(q)} 
\ll 1/\sqrt{\log \log q}$. 
Thus our Theorem gives a genuine asymptotic formula for all large $q$.

For any character $\chi\pmod q$ (not necessarily primitive) 
let $\frak a=0$ or $1$ be given by $\chi(-1)=(-1)^{\frak a}$. 
For $x >0$ we define
$$
W_{\frak a} (x) =  
\frac{1}{2\pi i} \int_{c-i\infty}^{c+i\infty} 
\Big(\frac{\Gamma(\tfrac {s+\frac 12 +\frak a}{2})}
{\Gamma(\tfrac {\frac 12+\frak a}{2})}\Big)^2 x^{-s} \frac{ds}{s}, 
\tag{1.2}
$$
for any positive $c$.  By moving the line of integration to $c=-\tfrac12 +\epsilon$ 
we may see that 
$$
W(x) = 1+ O(x^{\frac 12 -\epsilon}), \tag{1.3a} 
$$
and from the definition (1.2) we also get that 
$$
W(x) = O_c (x^{-c}). \tag{1.3b}
$$
We define 
$$
A(\chi) := \sum\Sb a, b =1\endSb^{\infty}  \frac{\chi(a)\cbar(b)}
{\sqrt{ab}} W_{\frak a}\Big(\frac{\pi ab}{q}\Big). \tag{1.4}
$$
If $\chi$ is primitive then $|L(\tfrac 12,\chi)|^2 = 2A(\chi)$ 
(see Lemma 2 below).  Let $Z=q/2^{\omega (q)}$ and decompose $A(\chi)$ as 
$B(\chi)+C(\chi)$ where 
$$
B(\chi)=\sum\Sb a, b \ge 1\\ ab\le Z \endSb \frac{\chi(a)\cbar(b)}
{\sqrt{ab}} W_{\frak a}\Big(\frac{\pi ab}{q}\Big),  
$$
and 
$$
C(\chi) = \sum\Sb a, b \ge 1\\ ab > Z \endSb \frac{\chi(a)\cbar(b)}
{\sqrt{ab}} W_{\frak a}\Big(\frac{\pi ab}{q}\Big). 
$$
Our main theorem will follow from the following two Propositions. 

\proclaim{Proposition 1} We have 
$$
\sumstar_{\chi \pmod q} |B(\chi)|^2 = \frac{\phi^*(q)}{8\pi^2} 
\prod_{p|q} \frac{(1-1/p)^3}{(1+1/p)} (\log q)^4 
\Big( 1+ O\Big(\frac{\omega(q)}{\log q}\Big) 
\Big).
$$
\endproclaim 

\proclaim{Proposition 2} We have 
$$
\sum_{\chi \pmod q} |C(\chi)|^2 \ll q \Big(\frac{\phi(q)}{q}\Big)^5 
(\omega(q) \log q)^2 + 
q (\log q)^3.
$$
\endproclaim

\demo{Proof of the Theorem} Since $|L(\tfrac 12,\chi)|^2 =2A(\chi)=2 (B(\chi)+C(\chi))$ 
for primitive characters $\chi$ we have 
$$
\sumstar_{\chi \pmod q} |L(\tfrac 12,\chi)|^4 
= 4 \sumstar_{\chi \pmod q} \Big(|B(\chi)|^2 + 2 B(\chi) C(\chi) + |C(\chi)|^2\Big).
$$
The first and third terms on the right hand side are handled 
directly by Propositions 1 and 2.  By Cauchy's inequality
$$
\sumstar_{\chi \pmod q} |B(\chi) C(\chi)| 
\le \Big(\sumstar_{\chi\pmod q} |B(\chi)|^2\Big)^{\frac 12} 
\Big(\sum_{\chi\pmod q} |C(\chi)|^2 \Big)^{\frac 12},
$$
and thus Propositions 1 and 2 furnish an estimate for the second term also.  
Combining these results gives the Theorem.

\enddemo

In [3], Heath-Brown refined Ingham's fourth moment for $\zeta(s)$, 
and obtained an asymptotic formula with a remainder term $O(T^{\frac 78+\epsilon})$. 
It remains a challenging open problem to obtain an asymptotic 
formula for $\sumstar_{\chi \pmod q} |L(\tfrac 12,\chi)|^4$ where 
the error term is $O(q^{1-\delta})$ for some positive $\delta$. 

This note arose from a conversation with Roger Heath-Brown at the 
Gauss-Dirichlet conference where he reminded me of this problem.  
It is a pleasure to thank him for this and other stimulating 
discussions.

\head 2. Lemmas \endhead

\proclaim{Lemma 1} If $(r,q)=1$ then 
$$
\sumstar_{\chi \pmod q} \chi(r)  
= \sum_{k|(q,r-1)} \phi(k) \mu(q/k).
$$
\endproclaim 
\demo{Proof} If we write $h_r(k) = \sumstar_{\chi \pmod k} \chi(r)$ 
then for $(r,q)=1$ we have 
$$
\sum_{k|q} h_r(k) = \sum_{\chi \pmod q} \chi(r) = \cases 
\phi(q) &\text{if  } q|r-1 \\ 
0 &\text{otherwise}.\\
\endcases 
$$
The Lemma now follows by M{\" o}bius inversion. 
\enddemo 

Note that taking $r=1$ gives the 
formula for $\phi^*(q)$ given in the introduction.  
If we restrict attention to characters of a given sign $\frak a$ 
then we have, for $(mn,q)=1$,
$$
\sumstar\Sb \chi \pmod q \\ \chi(-1)=(-1)^{\frak a}\endSb 
\chi(m)\cbar(n) 
= \frac 12 \sum_{k|(q,|m-n|)} \phi(k)\mu(q/k) + 
\frac{(-1)^{\frak a}}{2} \sum_{k|(q,m+n)} \phi(k) \mu(q/k). 
\tag{2.1}
$$

\proclaim{Lemma 2} If $\chi$ is a primitive character $\pmod q$ 
with $\chi(-1)=(-1)^{\frak a}$ then 
$$
|L(\tfrac 12,\chi)|^2  = 2 A(\chi),
$$
where $A(\chi)$ is defined in (1.4).
\endproclaim 

\demo{Proof} We recall the functional equation (see Chapter 9 of [1])
$$
\Lam(\tfrac 12 + s,\chi) = \Big(\frac {q}{\pi}\Big)^{s/2} \Gamma\Big(
\frac{s+\frac 12+ \frak a}{2}\Big) L(\tfrac 12+ s,\chi) 
= \frac{\tau(\chi)}{i^{\frak a}\sqrt{q}}
\Lam(\tfrac 12-s,\cbar),
$$
which yields 
$$
\Lam(\tfrac 12+s,\chi)\Lam(\tfrac 12+ s,\cbar)=\Lam(\tfrac 12 -s,\chi)\Lam(\tfrac 12-s,\cbar). \tag{2.2}
$$
For $c>\tfrac 12$ we consider 
$$
I:= \frac{1}{2\pi i} \int_{c-i\infty}^{c+i\infty} 
\frac{\Lam(\frac 12+s,\chi)\Lam(\frac 12+s,\cbar) }{\Gamma(\frac{\frac 12+\frak a}{2})^2} 
\frac{ds}{s}.
$$
We move the line of integration to Re$(s)=-c$, and use the 
functional equation (2.2).  This readily gives that 
$I=|L(\tfrac 12,\chi)|^2 -I$, so that $|L(\tfrac 12,\chi)|^2 =2I$.  
On the other hand, expanding $L(\tfrac 12+s,\chi)L(\tfrac 12+s,\cbar)$ into 
its Dirichlet series and integrating termwise, we get that 
$I=A(\chi)$.  This proves the Lemma.

\enddemo

We shall require the following bounds for divisor sums.  
If $k$ and $\ell$ are positive integers with 
$\ell k \ll x^{\frac 54}$ then
$$
\sum\Sb n\le x \\ (n,k)=1\endSb 
d(n) d(\ell k\pm n) \ll x(\log x)^2 \sum_{d|\ell} d^{-1}, \tag{2.3}
$$ 
provided that $x\le \ell k$ if the negative sign holds.  
This is given in (17) of Heath-Brown [2]. 
Secondly we record a result of P. Shiu [4] which 
gives that 
$$
\sum\Sb n\le x\\ n\equiv r\pmod k \endSb d(n) \ll 
\frac{\phi(k)}{k^2} x\log x, \tag{2.4}
$$
where $(r,k)=1$ and $x\ge k^{1+\delta}$ for some fixed $\delta >0$.

\proclaim{Lemma 3} Let $k$ be a positive integer, and let $Z_1$ and $Z_2$ 
be real numbers $\ge 2$.  If $Z_1Z_2 > k^{\frac {19}{10}}$ then 
$$
\sum\Sb Z_1 \le ab < 2Z_1\\ Z_2 \le cd < 2Z_2 \\ (abcd,k)=1\\ ac \equiv \pm bd \pmod k \\ 
ac\neq bd \endSb 1 \ll \frac{Z_1 Z_2}{k} (\log (Z_1 Z_2))^3.
$$
If $Z_1Z_2 \le k^{\frac {19}{10}}$ the quantity 
estimated above is $\ll (Z_1 Z_2)^{1+\epsilon}/k$.  
\endproclaim 
\demo{Proof} By symmetry we may just focus on the terms with $ac>bd$.  
Write $n=bd$ and $ac= k\ell \pm bd$.  Note that $k\ell \le 2ac$ and 
so $1\le \ell \le 8Z_1Z_2/k$.  Moreover since $ac\ge k\ell/2$ we 
have that $bd \le 4Z_1Z_2/(ac) \le 8Z_1Z_2/(k\ell)$.  Thus the 
sum we desire to estimate is 
$$
\ll \sum\Sb 1\le \ell \le 8Z_1Z_2/k \endSb \ \ \ \ \sum\Sb n\le 8Z_1Z_2/(k\ell) \\ n < k\ell\pm n \\
(n,k)=1\endSb 
d(n) d(k\ell \pm n). \tag{2.5}
$$
Since $d(n)d(k\ell \pm n) \ll (Z_1Z_2)^{\epsilon}$ the 
second assertion of the Lemma follows.  

Now suppose that $Z_1Z_2 > k^{\frac{19}{10}}$. 
We distinguish the cases $k\ell \le (Z_1Z_2)^{\frac {11}{20}}$ and 
$k\ell > (Z_1Z_2)^{\frac{11}{20}}$.  In the first case we estimate 
the sum over $n$ using (2.3).  Thus such terms contribute to (2.5) 
$$
\ll \sum_{\ell \le (Z_1Z_2)^{\frac{11}{20}}/k} 
\frac{Z_1Z_2}{k \ell} (\log Z_1Z_2)^2 \sum_{d|\ell} d^{-1} 
\ll \frac{Z_1Z_2}{k} (\log Z_1Z_2)^3.
$$
Now consider the second case.  Here we sum over $\ell$ first.  Writing 
$m=k\ell \pm n (=ac)$ we see that such terms contribute 
$$
\ll \sum\Sb n\le 8Z_1Z_2/k\endSb d(n) \sum\Sb (Z_1Z_2)^{\frac {11}{20}}/2 
\le m \le 4Z_1Z_2/n \\ m \equiv \pm n \pmod k \endSb d(m),  
$$
and by (2.4) (which applies as $(Z_1Z_2)^{\frac {11}{20}} > 
k^{\frac{209}{200}}$) this is 
$$
\ll \sum_{n \le 8Z_1Z_2/k } d(n) \frac{Z_1Z_2}{kn} \log Z_1 Z_2 
\ll \frac{Z_1 Z_2}{k} (\log Z_1Z_2)^3.
$$
The proof is complete.
 
\enddemo 

The next two Lemmas are standard; we have provided brief 
proofs for completeness.

\proclaim{Lemma 4} Let $q$ be a positive integer and $x\ge 2$ be a 
real number.  Then 
$$
\sum\Sb n\le x\\ (n,q)=1\endSb \frac{1}{n} = \frac{\phi(q)}{q} 
\Big( \log x+ \gamma +\sum_{p|q} \frac{\log p}{p-1}\Big) 
+ O\Big( \frac{2^{\omega(q)}\log x}{x}\Big).
$$
Further $\sum_{p|q} \log p/(p-1) \ll 1 +\log \omega(q)$.  
\endproclaim 
\demo{Proof}   We have 
$$
\align
\sum\Sb n\le x\\ (n,q)=1\endSb \frac{1}{n} &= \sum\Sb d| q \endSb \mu(d) 
\sum\Sb n\le x \\ d|n\endSb \frac{1}{n} 
= \sum\Sb d |q \\ d\le x\endSb \frac{\mu(d)}{d} \Big(\log \frac xd 
+\gamma +O\Big(\frac dx\Big)\Big)\\
&= \sum_{d|q} \frac{\mu(d)}{d} \Big( \log \frac xd +\gamma\Big) 
+ O\Big(\frac{2^{\omega(q)}\log x}{x}\Big).
\endalign
$$
Since $-\sum_{d|q} (\mu(d)/d)\log d = \phi(q)/q \sum_{p|q} (\log p)/(p-1)$ 
the first statement of the Lemma follows.  Since $\sum_{p|q} \log p/(p-1)$ 
is largest when the primes dividing $q$ are the first $\omega(q)$ primes, 
the second assertion of the Lemma holds. 

\enddemo 

\proclaim{Lemma 5} We have 
$$
\sum\Sb n\le q \\ (n,q)=1\endSb \frac{2^{\omega(n)}}{n} 
\ll \Big(\frac{\phi(q)}{q}\Big)^2 (\log q)^2.
$$
For $x\ge \sqrt{q}$ we have 
$$
\sum\Sb n\le x\\ (n,q)=1\endSb 
\frac{2^{\omega(n)}}{n} \Big(\log \frac xn\Big)^2 
= \frac{(\log x)^4}{12 \zeta(2)} \prod_{p|q} \Big(\frac{1-1/p}{1+1/p}\Big) 
\Big(1 + O\Big(\frac{1+\log \omega(q)}{\log q}\Big)\Big).
$$
\endproclaim
\demo{Proof}  Consider for Re$(s)>1$ 
$$
F(s)= \sum\Sb n=1\\ (n,q)=1\endSb^{\infty} 
\frac{2^{\omega(n)}}{n} 
= \frac{\zeta(s)^2}{\zeta(2s)} \prod_{p|q} \frac{1-p^{-s}}{1+p^{-s}}.
$$
Since 
$$
\sum\Sb n\le q \\ (n,q) =1 \endSb 
\frac{2^{\omega(n)}}{n} \le e\sum\Sb n=1 \\(n,q)=1\endSb^{\infty} 
\frac{2^{\omega(n)}}{n^{1+1/\log q}} = eF(1+1/\log q),
$$
the first statement of the Lemma follows. To prove the second statement
we note that, for $c>0$, 
$$
\sum\Sb n\le x\\ (n,q)=1\endSb 
\frac{2^{\omega(n)}}{n} \Big(\log\frac xn\Big)^2 
= \frac{2}{2\pi i} \int_{c-i\infty}^{c+i\infty} F(1+s) \frac{x^s}{s^3}
ds.
$$
We move the line of integration to $c=-\tfrac 12+\epsilon$ and 
obtain that the above is 
$$
2 \mathop{\text{Res}}_{s=0}  \ \ F(1+s)\frac{x^{s}}{s^3} 
+ O(x^{-\tfrac 12+\epsilon} q^{\epsilon}).
$$
A simple residue calculation then gives the Lemma. 
\enddemo 

\head 3. Proof of Proposition 1\endhead 

\noindent Applying (2.1) we easily obtain that 
$$
\sumstar_{\chi \pmod q} |B(\chi)|^2 
= M+E, 
$$
where
$$
M:= \frac{\phi^*(q)}{2} 
\sum\Sb a,b, c, d \ge 1\\ ab \le Z, cd \le Z \\ ac = bd \\ (abcd, q)=1\endSb 
\frac{1}{\sqrt{abcd}}
\Big( W_0\Big(\frac{\pi ab}{q}\Big) W_0\Big(\frac{\pi cd}{q} 
\Big) + W_1\Big(\frac{\pi ab}{q}\Big) W_1\Big(\frac{\pi cd}{q} 
\Big)\Big) \tag{3.1}
$$
and 
$$
E = \sum_{k|q} \phi(k) \mu^2 (q/k) E(k), 
$$
with 
$$
E(k) \ll \sum\Sb (abcd,q)=1\\ k| (ac\pm bd) \\ ac\neq bd \\ ab, cd \le Z
\endSb \frac{1}{\sqrt{abcd}}.
$$

To estimate $E(k)$ we divide the terms $ab$, $cd\le Z$ into dyadic 
blocks.  Consider the block $Z_1\le ab <2Z_1$, and $Z_2\le cd <2Z_2$.
By Lemma 3 the contribution of this block to $E(k)$ is, if $Z_1Z_2 
>k^{\frac{19}{10}}$,
$$
\ll \frac{1}{\sqrt{Z_1Z_2}} \frac{Z_1 Z_2}{k} (\log Z_1Z_2)^3 
\ll \frac{\sqrt{Z_1Z_2}}{k} (\log q)^3,
$$
and is $\ll (Z_1 Z_2)^{\frac 12+\epsilon}/k$ 
if $Z_1Z_2 \le k^{\frac{19}{10}}$.
Summing over all such dyadic blocks we obtain that 
$E(k) \ll (Z/k) (\log q)^3 + k^{-\frac{1}{20}+\epsilon}$, and so 
$$
E\ll Z 2^{\omega(q)} (\log q)^3 \ll q(\log q)^3.
$$

We now turn to the main term (3.1).  If $ac=bd$ then we may 
write $a=gr$, $b=gs$, $c=hs$, $d=hr$, where $r$ and $s$ are 
coprime.  We put $n=rs$, and note that given $n$ there are 
$2^{\omega(n)}$ ways of writing it as $rs$ with $r$ and $s$ coprime.
Note also that $ab=g^2 rs=g^2 n$, and $cd= h^2 rs =h^2 n$.  Thus 
the main term (3.1) may be written as 
$$
M= \frac{\phi^*(q)}{2} \sum_{\frak a=0, 1} 
\sum\Sb n \le Z\\ (n,q)=1 \endSb \frac{2^{\omega(n)}}{n} 
\Big(\sum\Sb g \le \sqrt{Z/n} \\ (g,q)=1\endSb \frac{1}{g} 
W_{\frak a} \Big(\frac{\pi g^2 n}{q}\Big) \Big)^2. 
$$
By (1.3a) we have that 
$W_{\frak a} (\pi g^2 n/Z) = 1+ O(\sqrt{g} n^{\frac 14}/q^{\frac 14})$, 
and using this above we see that 
$$
M= \phi^*(q) \sum\Sb n\le Z\\ (n,q)=1 \endSb \frac{2^{\omega(n)}}{n} 
\Big(\sum\Sb g\le \sqrt{Z/n} \\ (g,q) =1 \endSb \frac 1g + 
O(2^{-\omega(q)/4})\Big)^2. 
$$

We split the terms $n\le Z$ into the cases $n\le Z_0$ and $Z_0 < 
n\le Z$, where we set $Z_0 = Z/9^{\omega(q)} = q/18^{\omega(q)}$.  
In the first case, Lemma 4 gives that the sum over $g$ is 
$(\phi(q)/q) \log \sqrt{Z/n} +O(1+\log \omega(q))$.  Thus the 
contribution of such terms to $M$ is 
$$
\align
&\phi^*(q) \sum\Sb n\le Z_0 \\ (n,q) =1 \endSb \frac{2^{\omega(n)}}{n} 
\Big( \frac{\phi(q)}{2q} \log \frac{Z}{n} + O(1+ \log \omega(q))\Big)^2
\\
=
&\phi^*(q) \Big(\frac{\phi(q)}{2q}\Big)^2 
\sum\Sb n\le Z_0 \\ (n,q)=1\endSb \frac{2^{\omega(n)}}{n} 
\Big( \Big(\log \frac{Z_0}{n}\Big)^2 + O(\omega(q) \log q)\Big).
\\
\endalign
$$
Using Lemma 5 we conclude that the terms $n\le Z_0$ contribute 
to $M$ an amount
$$
\frac{\phi^*(q)}{8\pi^2} 
 \prod_{p|q} \frac{(1-1/p)^3}{(1+1/p)}  (\log q)^4 \Big( 1 +O\Big(\frac{\omega(q)}{\log q}\Big)\Big).
 \tag{3.2}
$$  

In the second case when $Z_0 \le n \le Z$, 
we extend the sum over $g$ to all $g\le 3^{\omega(q)}$ that are 
coprime to $q$, and so by Lemma 4 the sum over 
$g$ is $\ll \omega(q)\phi(q)/q$.  Thus 
these terms contribute to $M$ an amount 
$$
\ll \phi^*(q) \Big(\omega(q) \frac{\phi(q)}{q}\Big)^2 \sum\Sb Z_0\le 
n\le Z\endSb \frac{2^{\omega(n)}}{n} \ll 
\phi^*(q) \Big(\frac{\phi(q)}{q}\Big)^2 (\omega(q))^3 \log q.
$$
Since $q\omega(q)/\phi(q) \ll \log q$, combining this with (3.2) 
we conclude that 
$$
M = \frac{\phi^*(q)}{8\pi^2} 
 \prod_{p|q} \frac{(1-1/p)^3}{(1+1/p)}  (\log q)^4 
\Big( 1 +O\Big(\frac{\omega(q)}{\log q}\Big)\Big).
$$
Together with our bound for $E$, this proves Proposition 1. 

\head 4. Proof of Proposition 2 \endhead

\noindent The orthogonality relation for characters gives that 
$$
\align
\sum_{\chi \pmod q} |C(\chi)|^2 
&\ll \phi(q) \sum\Sb (abcd,q)= 1\\ ac\equiv \pm bd \pmod q\\ ab, cd >Z \endSb 
\frac{1}{\sqrt{abcd}} \sum_{\frak a =0, 1} \Big|W_{\frak a} \Big( \frac{\pi ab}{q}\Big) W_{\frak a} 
\Big(\frac{\pi cd}{q}\Big)\Big|\\
&\ll \phi(q) \sum\Sb (abcd,q)= 1\\ ac\equiv \pm bd \pmod q\\ ab, cd >Z \endSb 
\frac{1}{\sqrt{abcd}} \Big( 1+ \frac{ab}{q}\Big)^{-2} \Big( 1+ \frac{cd}{q}\Big)^{-2},\\
\endalign
$$
using (1.3a,b).  We write the last expression above as $R_1 + R_2$ where $R_1$ 
contains the terms with $ac=bd$, and $R_2$ contains the rest.  

We first get an estimate for $R_2$.  We break up the terms into dyadic blocks; a 
typical one counts $Z_1 \le ab < 2Z_1$ and $Z_2 \le cd < 2Z_2$ (both $Z_1$ and $Z_2$ 
being larger than $Z$).  The contribution of such a dyadic block is, using 
Lemma 3, (note that $Z_1Z_2 >Z^2 >q^{\frac{19}{10}}$)
$$
\ll \frac{\phi(q)}{\sqrt{Z_1Z_2}} \Big(1+ \frac{Z_1}{q}\Big)^{-2} 
\Big( 1+ \frac{Z_2}{q}\Big)^{-2} \frac{Z_1Z_2}{q} (\log Z_1Z_2)^3.
$$
Summing this estimate over all the dyadic blocks we obtain that 
$$
R_2 \ll q (\log q)^3. 
$$

We now turn to the terms $ac=bd$ counted in $R_1$.  As in our 
treatment of $M$, we write $a=gr$, $b=gs$, $c=hs$, $d=hr$, with $(r,s)=1$, 
and group terms according to $n=rs$.  We see easily that 
$$
R_1 \ll \phi(q) \sum\Sb (n,q) =1\endSb \frac{2^{\omega(n)}}{n}  
\Big( \sum\Sb g >\sqrt{Z/n} \\ (g,q) =1 \endSb  
\frac{1}{g} \Big(1+ \frac{g^2 n}{q}\Big)^{-2} \Big)^2. \tag{4.1} 
$$
First consider the terms $n> q$ in (4.1).  Here the sum over $g$ gives an 
amount $\ll q^{2}/n^2$ and so the contribution of these terms 
to (4.1) is 
$$
\ll \phi(q) \sum\Sb n> q  \endSb \frac{2^{\omega(n)}}{n} 
\frac{q^4}{n^4} 
\ll \phi(q) \log q.
$$
For the terms $n<q$ the sum over $g$ in (4.1) is 
easily seen to be 
$$
\ll 1+ \sum\Sb \sqrt{Z/n} \le g\le \sqrt{q/n} \\ (g,q) =1 \endSb 
\frac{1}{g} \ll 1 + \frac{\phi(q)}{q} \omega(q).
$$
The last estimate follows from Lemma 4 when $n < Z/9^{\omega(q)}$, 
while if $n>Z/9^{\omega(q)}$ we extend the sum over $g$ to 
all $g\le 6^{\omega(q)}$ with $(g,q)=1$ and then use Lemma 4. 
Thus the contribution of terms $n<q$ to (4.1) is, using Lemma 5, 
$$
\ll \phi(q) \Big(1+\frac{\phi(q)}{q}\omega(q)\Big)^2 
\sum\Sb n\le q \\ (n,q) =1\endSb \frac{2^{\omega(n)}}{n} 
\ll q \log^2 q \Big(\frac{\phi(q)}{q}\Big)^5 \omega(q)^2.
$$
Combining these bounds with our estimate for $R_2$ we obtain Proposition 2.

\enddocument